\documentclass [leqno, 10pt] {amsart}

\usepackage[all]{xy}
\usepackage{amssymb, amsmath, amsfonts, amsthm, amsbsy, amscd}

\newtheorem{theorem}{Theorem}[section]
\newtheorem{proposition}{Proposition}[section]
\newtheorem{corollary}{Corollary}[section]
\newtheorem{lemma}{Lemma}[section]

\newtheorem{definition}{Definition}[section]

\theoremstyle{definition}

\newtheorem{remark}{Remark}[section]

\numberwithin{equation}{section}

\newcommand{\Ct}{\mathcal{C}}
\newcommand{\St}{\mathcal{S}}
\newcommand{\Tt}{\mathcal{T}}
\newcommand{\bun}{\mathcal{F}}
\newcommand{\emfz}{\mathcal{E}^{\circ}}
\newcommand{\stfree}{\mathcal{E}^{\circ}_{(ij)}\,^{k}}
\newcommand{\st}{\mathcal{E}_{(ij)}\,^{k}}

\newcommand{\dn}{d^{\eta}}
\newcommand{\hschw}{\Hat{\mathbf{S}}}
\newcommand{\hschwn}{\Hat{\mathbf{S}}_{[\nabla]}\,}
\newcommand{\hpnabla}{\widehat{^\phi\nabla}}
\newcommand{\nablap}{\nabla^{\prime}}

\newcommand{\pL}{^{\phi}L}
\renewcommand{\xi}{\frac{\partial}{\partial x^{i}}}

\newcommand{\pnabla}{^{\phi}\nabla}

\newcommand{\schw}{\mathbf{S}}
\newcommand{\schwn}{\mathbf{S}_{[\nabla]}\,}

\newcommand{\ab}{\mathcal{A}}

\newcommand{\standrepz}{\mathbb{V}^{\times}}

\newcommand{\contacto}{\mathcal{CO}}

\newcommand{\rep}{\mathbb{T}}

\newcommand{\G}{\mathcal{G}}

\newcommand{\emf}{\mathcal{E}}

\newcommand{\hnabla}{\hat{\nabla}}

\newcommand{\standrep}{\mathbb{V}}

\newcommand{\spl}{\mathfrak{sp}}

\newcommand{\form}{\mathsf{L}}

\newcommand{\eul}{\mathbb{X}}

\newcommand{\rb}{T}
\newcommand{\proj}{\mathbb{P}}

\newcommand{\lie}{\mathfrak{L}}

\newcommand{\vect}{\text{Vec}}

\newcommand{\g}{\mathfrak{g}}

\newcommand{\p}{\mathfrak{p}}

\newcommand{\tensor}{\otimes}

\newcommand{\rea}{\mathbb R}

\begin{document}
\title{Contact Schwarzian Derivatives}
\author{Daniel J. F. Fox} 
\address{School of Mathematics\\ Georgia Institute of Technology\\ 686 Cherry St.\\ Atlanta, GA 30332-0160, U.S.A.}
\email{fox@math.gatech.edu}
\begin{abstract}
H. Sato introduced a Schwarzian derivative of a contactomorphism of $\rea^{3}$ and with T. Ozawa described its basic properties. In this note their construction is extended to all odd dimensions and to non-flat contact projective structures. The contact projective Schwarzian derivative of a contact projective structure is defined to be a cocycle of the contactomorphism group taking values in the space of sections of a certain vector bundle associated to the contact structure, and measuring the extent to which a contactomorphism fails to be an automorphism of the contact projective structure. For the flat model contact projective structure, this gives a contact Schwarzian derivative associating to a contactomorphism of $\rea^{2n-1}$ a tensor which vanishes if and only if the given contactomorphism is an element of the linear symplectic group acting by linear fractional transformation.
\end{abstract}

\maketitle
\section{Introduction}
The classical Schwarzian derivative is a cocycle, $\schw(f)$, of the diffeomorphism group of the real line with coefficients in the quadratic differentials and which vanishes when restricted to the group of projective transformations. It arises naturally in the context of flat projective structures on one-dimensional manifolds and may be interpreted as describing how a normalized second order linear differential operator transforms under a change of variable, provided the operator is viewed as acting on $-1/2$ densities, rather than on functions (this determines the relevant normalization). Its characteristic properties are:
\begin{enumerate}
\item $\schw(f) = 0$ if and only if $f$ is the restriction of a linear fractional transformation.
\item $\schw(f \circ g) = g^{\ast}(\schw(f)) + \schw(g)$. (Cocycle property).
\item Locally there are linearly independent solutions, $x_{1}$ and $x_{2}$, of $\ddot{x} + \tfrac{1}{2}\schw(f)x = 0$, such that $f = \frac{x_{2}}{x_{1}}$.\end{enumerate}

This note describes the generalized Schwarzian derivatives associated to contactomorphisms by the action of the contactomorphism group on the space of contact projective structures on a contact manifold. (See \cite{Fox} for background on contact projective structures). The basic properties of these generalized Schwarzian derivatives directly generalize those of the classical Schwarzian derivative. In \cite{Sato}, H. Sato introduced a Schwarzian derivative of a contactomorphism of $\rea^{3}$ and with T. Ozawa, in \cite{Ozawa-Sato}, he explored its basic properties. Their Schwarzian derivative measures the failure of a contactomorphism to preserve the flat model contact projective structure in three dimensions. In this note their construction is extended to all odd dimensions and to non-flat contact projective structures. The contact projective Schwarzian derivative of a contact projective structure is defined to be a cocycle of the contactomorphism group taking values in the space of sections of a certain vector bundle associated to the contact structure, and measuring the extent to which a contactomorphism fails to be an automorphism of the contact projective structure. For the flat model contact projective structure, this gives a contact Schwarzian derivative associating to a contactomorphism of $\rea^{2n-1}$ a tensor which vanishes if and only if the given contactomorphism is an element of the linear symplectic group acting by linear fractional transformation. Even in the three-dimensional case, this point of view simplifies the construction of the contact Schwarzian and makes the proofs of its basic properties simpler than the corresponding proofs in \cite{Ozawa-Sato}.

The Schwarzian cocycles associated to different contact projective structures on the same contact manifold are cohomologous, and determine a canonical non-trivial class in the first cohomology of the contactomorphism group with coefficients in the space of sections of a certain vector bundle. This cohomology class depends only on the contact structure. In Section \ref{contactpathsection} the contact Schwarzians are discussed briefly in the more general setting of contact path geometries, which is closer to the original point of view of Sato. The main conclusion is that a contactomorphism is completely determined locally by a single function on the total space of the projectivized contact hyperplane bundle and satisfying some complicated integrabiliy condition. Many of the results presented here have parallels or specializations in Ozawa-Sato, \cite{Ozawa-Sato}; as the translation is in general straightforward, though sometimes computationally involved, it has been left in general to the reader.

Generalized Schwarzian derivatives have been studied by various authors in the contexts of projective structures and conformal structures. The approach to contact Schwarzian derivatives taken here is basically a generalization of M. Yoshida's point of view on projective Schwarzian derivatives (see the survey, \cite{Sasaki-Yoshida}, of T. Sasaki - Yoshida for applications and more complete references), though informed by the general theory of parabolic geometries. The essence of some aspects of this construction was understood already by T. Y. Thomas, \cite{Thomas}, in the 1920's, and subsequently by those studying the invariant differential operators arising in the context of parabolic geometries, see e.g. \cite{Baston-Eastwood}. The idea of regarding a Schwarzian derivative as a cocycle of some group of diffeomorphisms goes back at least to R. C. Gunning, e.g. \cite{Gunning-3} and \cite{Gunning}, and has figured prominently recently in the papers of S. Bouarroudj, C. Duval, C. Lecomte, and V. Ovsienko, e.g. \cite{Bouarroudj-Ovsienko}, \cite{Duval-Ovsienko}, \cite{Lecomte-Ovsienko}, \cite{Lecomte-Ovsienko-modules}. See also the forthcoming textbook of Ovsienko and S. Tabachnikov, \cite{Ovsienko-Tabachnikov}, for a detailed presentation of the projective Schwarzian derivative. There are many other sources, but no effort has been made to survey them here.

\section{Background}
\subsection{Preliminaries and Notations}\label{notationsection}
Let $(M, H)$ be a $(2n-1)$-dimensional contact manifold. Each choice, $\theta$, of a contact one-form determines uniquely a Reeb vector field characterized by $\theta(\rb) = 1$ and $i(\rb)d\theta = 0$. The choice of contact one-form is refered to as a choice of \textbf{scale}. Lowercase Latin indices will run from $1$ to $2n-2$. Lowercase Greek indices will run from $0$ to $2n-2$. A coframe, $\theta^{\alpha}$, is \textbf{$\theta$-adapted} if $\theta^{0} = \theta$ and $\theta^{i}(\rb) = 0$. An adapted coframe determines a dual frame, $E_{\alpha}$, such that $E_{0} = \rb$ and the $E_{i}$ span $H$. When a contact form is fixed, an adapted coframe and corresponding dual frame will be assumed fixed also. The notations $S_{[\alpha_{1}\dots \alpha_{k}]}$ and $S_{(\alpha_{1}\dots \alpha_{k})}$ denote, respectively, the complete skew-symmetrization and the complete symmetrization over the bracketed indices. Sometimes the abstract index notation will be used, so that equations with indices have invariant meaning. Greek abstract indices label sections of tensor bundles on $M$, while Latin abstract indices label sections of the tensor powers of $H$ and $H^{\ast}$, so that an expression such as $\tau_{[ij]}\,^{k}$ indicates a section of $\Lambda^{2}(H^{\ast})\tensor H$. Each $\theta$ determines a splitting, $TM = H \oplus \text{span}\{\rb\}$, which induces a splitting of the full tensor bundle. Using these splittings Latin abstract indices may be interpreted as the components of a tensor with respect to a $\theta$-adapted coframe and dual frame. The components of $\omega = d\theta$ are $\omega_{\alpha\beta} =  \omega_{[\alpha\beta]} = \omega(E_{\alpha}, E_{\beta})$. As $\omega_{0\alpha} = 0$, $\omega$ may be written as $\omega = \frac{1}{2}\omega_{ij}\theta^{i}\wedge \theta^{j}$. Latin indices may be raised and lowered using $\omega_{ij}$ according to the following conventions. Defining $\omega^{kl}$ by $\omega^{kl}\omega_{lj} = -\delta_{j}\,^{k}$, let $\gamma^{p} = \omega^{pq}\gamma_{q}$, and $\gamma_{p} = \gamma^{q}\omega_{qp}$. It is necessary to pay attention to which index is raised or lowered as, for instance, $\eta^{p}\gamma_{p} = -\eta_{p}\gamma^{p}$. Under a change of scale,  $\tilde{\theta} = f^{2}\theta$, ($f \neq 0$), the restriction to $H$ of $\omega$ rescales by $f^{2}$, so there is induced on $H$ a well-defined conformal symplectic structure. 
Expressions labeled with a $\tilde{\,\,}$ indicate use of a $\Tilde{\theta}$-adapted coframe and dual frame, unless indicated otherwise. 

In general $(M, H)$ will be assumed co-oriented, so that the bundle, $(TM/H)^{\ast}$, of contact one-forms has structure group $\rea^{>0}$, and $\form$ will denote the principal $\rea^{\times}$ bundle of frames in a chosen square-root of $(TM/H)^{\ast}$. A contact one-form consistent with the chosen co-orientation is called positive. Because under rescaling the contact volume transforms by $\tilde{\theta} \wedge (d\tilde{\theta})^{n-1} = f^{2n}\theta \wedge (d\theta)^{n-1}$, the bundle $\form$ is naturally identified with a $1/2n$th root of the bundle of frames in the canonical bundle $\wedge^{2n-1}(T^{\ast}M)$ having the orientation induced by the volume form associated to a positive contact one-form. Denote by $\emf[\lambda]$ the line bundle associated to $\form$ by the representation, $r\cdot s = r^{-\lambda}s$, of $\rea^{\times}$ on $\rea$, so that $\emf[-1]$ is a $1/2n$th root of $\Lambda^{2n-1}(T^{\ast}M)$. The model for $\form$ is the defining bundle $\standrepz \to \proj(\standrep)$, where $(\standrep, \Omega)$ is a real symplectic vector space. Notation such as $\emf_{(ij)}\,^{k}[\lambda]$ indicates the tensor product $S^{2}(H^{\ast})\tensor H \tensor \emf[\lambda]$, and the addition of a superscript, $\,^{\circ}\,$, indicates the subbundle $\stfree \subset \st$ comprising completely trace free sections. Let $\St$, $\Tt$, and $\Ct$, respectively, denote the bundles of tensors on $H$ obtained by raising the third index of elements of, respectively, $S^{3}(H^{\ast})$; the subbundle of $\tensor^{3}(H^{\ast})$ comprising trace free tensors satisfying $T_{i(jk)} = T_{ijk}$ and $T_{(ijk)} = 0$; and the subbundle of $\tensor^{3}(H^{\ast})$ comprising trace free tensors satisfying $C_{i(jk)} = C_{ijk}$. Though the operation of raising an index depends on the choice of contact one-form, the bundles so defined do not. By results of Weyl, \cite{Weyl}, the fiber over a point of any of $\St$, $\Tt$, and $\Ct$, is a semisimple $Sp(n-1, \rea)$-module, and $\Ct = \St \oplus \Tt$ is a decomposition into irreducibles. The notation $\Gamma(\St)$ denotes the space of smooth sections of $\St$.

\subsection{Statement of Main Theorem}
For the definitions of the structures involved in Theorem \ref{maintheorem} see Section \ref{definitionsection} below. 
\begin{theorem}\label{maintheorem}
Given a contact projective structure, $(M, H, [\nabla])$, there is associated to each contactomorphism, $\phi$, of $(M, H)$, a section, $\schwn(\phi) \in \Gamma(\Ct)$, having the following properties: \\
\noindent
$(1)$ $\schwn(\phi) = 0$ if and only if $\phi$ is an automorphism of $(M, H, [\nabla])$.\\
\noindent
$(2)$ For any contactomorphisms, $\phi$ and $\psi$, $\schwn(\phi \circ \psi) = \psi^{\ast}(\schwn(\phi)) + \schwn(\psi)$.\\
\noindent
$(3)$ There is a contact projectively invariant differential operator, $L:\emf[1] \to \emf_{(ij)}[1]$, so that if $\pL = \phi^{\ast} \circ L \circ (\phi^{-1})^{\ast}$, then
\begin{align*}
&\pL_{ij}u - L_{ij}u = -\schwn_{(ij)}\,^{k}\nabla_{k}u + \tfrac{1}{2n(3-2n)}\left(\nabla^{p}\schwn_{pij} - \schwn_{pqi}\schwn^{pq}\,_{j}\right)u\\&  + \tfrac{n-1}{n(2n-3)}\left(\nabla_{p}\schwn_{(ij)}\,^{p} - \schwn_{q(j}\,^{p}\schwn_{i)p}\,^{q} -\schwn_{q(j}\,^{p}\tau_{i)p}\,^{q}\right)u,
\end{align*}
where $\nabla \in [\nabla]$ is the unique representative of the given contact projective structure associated to a chosen contact one-form (as in Theorem \ref{theorema} below), and $\tau_{ij}\,^{k}$ is the contact torsion of $[\nabla]$. If the given contact projective structure has vanishing contact torsion, then $\schwn(\phi)\in\Gamma(\St)$, and the operator $L:\emf[1] \to \emf_{(ij)}[1]$, transforms as
\begin{align*}
&\pL_{ij}u - L_{ij}u = -\schwn_{ij}\,^{k}\nabla_{k}u  +\tfrac{1}{2n}(\nabla_{p}\schwn_{ij}\,^{p} - \schwn_{qj}\,^{p}\schwn_{ip}\,^{q})u.
\end{align*}
Moreover, the given contact projective structure is flat if and only if in a neighborhood of every point of $M$ the equation $L(u) = 0$ admits $2n$ linearly independent solutions. On a manifold with flat contact projective structure the equation $\pL(u) = 0$ admits locally $2n$ linearly independent solutions, $u^{1}, \dots, u^{2n}$, from which the contactomorphism $\phi$ may be reconstructed locally.
\end{theorem}

\noindent
This functional $\schwn(\phi)$ is the desired contact generalization of the classical Schwarzian derivative, and the remainder of this note is devoted to motivating its construction and proving its basic properties. The first and third claims of Theorem \ref{maintheorem} have the following specializations in the case that the given contact projective structure is the flat model:\\
\noindent
$(1^{\prime})$ $\schw(\phi)$ vanishes on an open set $U$ if and only if $\phi$ equals the restriction to $U$ of an element of the linear symplectic group $Sp(n, \rea)$ acting on $\proj(\rea^{2n})$.\\
\noindent
$(3^{\prime})$ 
Given any $\theta$-adapted coframe, $\theta^{\alpha}$, and dual frame, $E_{\alpha}$, parallel with respect to the unique representative, $\nabla \in [\nabla]$ associated to $\theta$ (as in Theorem \ref{theorema} below), there are locally $2n$ linearly independent sections of $\emf[1]$, $f^{\infty}, f^{1}, \dots, f^{2n-1}, f^{0}$, solving the system of PDE:
\begin{align}\label{schwarzianpde}
(E_{i}(E_{j}(f)) + E_{j}(E_{i}(f))) = 2\schw_{ij}\,^{p}E_{p}(f)  - \tfrac{1}{n}(E_{p}(\schw_{ij}\,^{p}) -\schw_{pi}\,^{q}\schw_{qj}\,^{p}),
\end{align}
and such that $\phi^{\alpha} = \frac{f^{\alpha}}{f^{\infty}}$.

\subsection{Review of Contact Projective Structures}\label{definitionsection}
In this section are summarized the needed facts about contact projective structures; proofs and further details may be found in \cite{Fox}. Call a smoothly immersed one-dimensional submanifold a {\bf path}. Call a path everywhere tangent to $H$ a {\bf contact path}. An affine connection, $\nabla$, is said to admit a \textbf{full set of contact geodesics} if every geodesic of $\nabla$ tangent to $H$ at one point is everywhere tangent to $H$. It is easily checked that $\nabla$ admits a full set of contact geodesics if and only if $\nabla_{(i}\theta_{j)} = 0$ for any choice of contact one-form, $\theta$.
\begin{definition}\label{imprecisedefncontactpath}
A {\bf contact path geometry} is a $(4n-5)$ parameter family of contact paths in $(M, H)$ such that for every $x \in M$ and each $L \in \proj(H_{x})$ there is in the family a unique contact path containing $x$ and tangent to $L$. Two such families of paths are {\bf equivalent} if there is a contactomorphism mapping the paths of one family onto the paths of the other family. A \textbf{contact projective structure} is a contact path geometry the contact paths of which are among the unparameterized geodesics of some affine connection; in this case the contact paths are called \textbf{contact geodesics}. Two contact projective structures are {\bf equivalent} if and only if they are equivalent as contact path geometries. 
\end{definition}
\noindent
The model contact projective structure is the family of contact lines comprising the images in the projectivization of a symplectic vector space of the two-dimensional isotropic subspaces. A contact projective structure is \textbf{flat} if it is locally equivalent to this model. 

\begin{theorem}[\cite{Fox}]\label{theorema}
Given a contact projective structure, there is associated to each choice of contact one-form, $\theta$, a unique affine connection, $\nabla$, with torsion tensor, $\tau$, having among its unparameterized geodesics the given contact geodesics and satisfying $\nabla\theta = 0$; $\nabla d\theta = 0$;  $\tau_{0i}\,^{\alpha} = 0 = \tau_{i0}\,^{\alpha}$; and $\tau_{ip}\,^{p} = 0 = \tau_{pi}\,^{p}$.
\end{theorem}
\begin{lemma}[\cite{Fox}]\label{transformationruleslemma}
Given a contact projective structure, let $\Lambda$ be the difference tensor of the representatives, $\tilde{\nabla}$ and $\nabla$, associated by Theorem \ref{theorema} to the choices of contact one-forms, $\tilde{\theta} = f^{2}\theta$ and $\theta$. With respect to a $\theta$-adapted coframe and dual frame, the components of $\Lambda$ are expressible in terms of $\gamma = d\log{f} = f^{-1}df$ as
\begin{align}
&\label{lambdatransform}\Lambda_{ij}\,^{k} = \gamma_{i}\delta_{j}\,^{k} + \gamma_{j}\delta_{i}\,^{k} + \omega_{ij}\gamma^{k},&
&\Lambda_{\alpha \beta }\,^{0} = 2\gamma_{\alpha}\delta_{\beta}\,^{0},\\
&\label{lambdai0p}\Lambda_{\alpha 0}\,^{j} = 4\gamma_{\alpha}\gamma^{j} - 2\nabla_{\alpha}\gamma^{j} + 4\delta_{\alpha}\,^{0}\gamma^{q}\nabla_{q}\gamma^{j},&
&\Lambda_{0i}\,^{j} = -2\gamma^{q}\tau_{qi}\,^{j} - 2\nabla_{i}\gamma^{j}.
\end{align}
\end{lemma}
\noindent
Given $\nabla$ as in Theorem \ref{theorema}, the components, $\tau_{ij}\,^{k}$, of the torsion of $\nabla$ do not depend the choice of scale, and this \textbf{contact torsion} is the most basic invariant of the contact projective structure. Let $R_{\alpha\beta\gamma}\,^{\sigma}$ denote the curvature of $\nabla$, note that $R_{\alpha\beta\gamma}\,^{0} = 0$, and let $R_{ij} = R_{ipj}\,^{p}$ and $S_{ij} = R_{p}\,^{p}\,_{ij}$ be the two possibly independent traces of its curvature tensor. By the contracted first Bianchi identity, $S_{ij} + 2R_{ij} = 2\nabla_{p}\tau^{p}\,_{ij} - \tau^{pq}\,_{j}\tau_{pqi}$, 
so $S_{ij} = -2R_{ij}$ if the contact torsion vanishes. The following tensors are basic in the study of contact projective structures:
\begin{align*}
&P_{ij} = \tfrac{1}{n(2n-3)}\left((n-1)R_{ij} - \tfrac{1}{2n-1}R_{[ij]} + \tfrac{1}{4}S_{ij}\right),\\
&Q_{ij} = \tfrac{1}{3-2n}\left(2R_{ij} + S_{ij} -\tfrac{4}{2n-1}R_{[ij]}\right),\\
&W_{ijk}\,^{l} = R_{ijk}\,^{l} + 2\delta_{[i}\,^{l}P_{j]k} +2\omega_{k[j}P_{i]}\,^{l} + 2\omega_{ij}P_{k}\,^{l} + \omega_{ij}Q_{k}\,^{l},\\
&C_{ijk} = R_{0ijk} - \left(2\nabla_{i}P_{jk} + \nabla_{i}Q_{jk}\right) - \tfrac{2}{2n-1}(2\omega_{i(k}\nabla^{p}P_{j)p} + \omega_{i(k}\nabla^{p}Q_{j)p}).
\end{align*}
\noindent
$W_{ijk}\,^{l}$ is the {\bf contact projective Weyl tensor}, and $C_{ijk}$ should be regarded as an analogue of the Cotton tensor in conformal geometry. There hold the following identities.
\begin{align*}
 &2(1-n)Q_{ij} + Q_{ji} = 2R_{ij} + S_{ij},& &Q_{ij} = 2R_{ij} - 4nP_{ij},& \\ &Q_{[ij]} = -2P_{[ij]} = -\tfrac{2}{2n-1}R_{[ij]}.&
\end{align*}
$R_{p}\,^{p} = 0$ and so $P_{p}\,^{p} = 0 = Q_{p}\,^{p}$. When the contact torsion vanishes, $P_{ij} = P_{(ij)} = \frac{1}{2n}R_{ij}$, $Q_{ij} = 0$, and $C_{ijk} = R_{0ijk} - 2\nabla_{i}P_{jk} - \tfrac{4}{2n-1}\omega_{i(j}\nabla^{p}P_{k)p}$. The definitions and the Bianchi identities imply the following identities.
\begin{align*}
&W_{p}\,^{p}\,_{ij} = 0, \quad &W_{ijp}\,^{p} = 0&, \quad &W_{ipj}\,^{p} = - \tfrac{1}{2}Q_{ij},&\\
&C_{ip}\,^{p} = 0, \quad &C_{p}\,^{p}\,_{k} = 0&, \quad &C_{i[jk]} = 0.
\end{align*}
$W_{ijk}\,^{l}$ is invariant if the contact torsion vanishes. Because $\tau_{[ijk]} = 0$ and $W_{[ijk]l} = 0$, basic facts (see \cite{Weyl}) about representations of $Sp(2, \rea)$ imply that in dimension three $\tau_{ij}\,^{k}$ and $W_{ijk}\,^{l}$ vanish identically. In three dimensions $C_{ijk}$ is invariant and completely symmetric. A contact projective structure is locally flat in dimension three if and only if $C_{ijk} = 0$, and in dimensions greater than three if and only if both $\tau_{ij}\,^{k} = 0$ and $W_{ijk}\,^{l} = 0$.

In Theorem \ref{ambienttheorem} the `ambient' manifold, $\rho:\form \to M$, is a square-root of the bundle of positive contact one-forms on the co-oriented contact manifold, $(M, H)$; $\eul$ is the vertical vector field generating the dilations in the fibers of $\form \to M$; $\alpha$ is the tautological one-form on $\form$ defined by $\alpha_{p}(X) = p^{2}(\rho_{\ast}(X))$; and $\Omega = d\alpha$ is the canonical symplectic structure on $\form$. Any (local) section, $s:M \to \form$, determines a horizontal lift, $\hat{X} \in \Gamma(T\form)$, of a vector field $X \in \Gamma(TM)$. Let uppercase Latin indices run over $\{\infty, 1, \dots, 2n-2, 0\}$, $\infty$ indicating the vertical direction. Denote by $\hat{R}_{IJK}\,^{L}$ the curvature tensor of an affine connection, $\hat{\nabla}$, on $\form$, and raise and lower indices with $\Omega_{IJ}$.
\begin{theorem}[\cite{Fox}]\label{ambienttheorem}
Let $(M, H)$ be a co-oriented contact manifold and let $\rho:\form \to M$ be a square-root of the bundle of positive contact one-forms. There is a functor associating to each contact projective structure on $M$ a unique affine connection, $\hat{\nabla}$, (the {\bf ambient connection}), on the total space of $\form$, having torsion $\hat{\tau}$, and satisfying: $(1)$. $\hat{\nabla}\eul$ is the fundamental $\binom{1}{1}$-tensor on $\form$; $(2)$. $i(\eul)\hat{\tau} = 0$; $(3)$. $\hat{\nabla}\Omega = 0$; $(4)$. The Ricci trace, $\hat{R}_{IPJ}\,^{P}$, of the curvature tensor of $\hat{\nabla}$ vanishes; $(5)$. There vanishes the restriction to $\ker \alpha$ of the tensor $\hat{R}_{Q}\,^{Q}\,_{IJ}$; $(6)$. The projections into $M$ of the unparametrized geodesics of $\hat{\nabla}$ transverse to the vertical and tangent to $\ker \alpha$ are the contact geodesics of the given contact projective structure. Moreover, the contact projective structures with vanishing contact torsion are in bijection with the torsion free affine connections satisfying conditions $(1)$, $(3)$, $(4)$, and $(6)$. In this case condition $(2)$ is vacuous, $(5)$ follows from $(4)$ by the contracted first Bianchi identity, and the curvature tensor is completely trace free. 
\end{theorem}
\noindent
Condition $(6)$ is equivalent to the following statement useful in computations: For any (local) section, $s:M \to \form$, the affine connection, $\bar{\nabla}$, on $M$ defined by $\bar{\nabla}_{X}Y = \rho_{\ast}(\hat{\nabla}_{\hat{X}}\hat{Y})$ represents the given contact projective structure.

\subsection{Flat Model Contact Projective Structure}\label{flatmodelsection}
Equip $\rea^{2n}$ with the symplectic form $\Omega = du^{\infty} \wedge du^{0} + \tfrac{1}{2}\omega_{pq}du^{p}\wedge du^{q} = \tfrac{1}{2}\Omega_{IJ}du^{I}\wedge du^{J}$, where again uppercase Latin indices run over $\{\infty, 1, \dots, 2n-2, 0\}$. Represent the general element of $G = Sp(n, \rea)$, as a matrix $A_{I}\,^{J}$ satisfying $A_{I}\,^{P}A_{J}\,^{Q}\Omega_{PQ} = \Omega_{IJ}$. $\proj(\rea^{2n})$ is a homogeneous space $G/P$, where $P$ is the stabilizer of a point in $\proj(\rea^{2n})$. In the chart on which $u^{\infty} \neq  0$, define coordinates by $x^{\alpha} = \frac{u^{\alpha}}{u^{\infty}}$. In these coordinates $G$ acts on $\proj(\rea^{2n})$ by linear fractional transformations, $x^{\alpha} \to \frac{A_{\alpha}\,^{\infty} + A_{\beta}\,^{\alpha}x^{\beta}}{A_{\infty}\,^{\infty} + A_{\beta}\,^{\infty}x^{\beta}}$. The right action of $G$ on $G$ induces a Lie algebra embedding, $\g = \spl(n, \rea) \to \vect(G)$, defined by $h \to X_{h}(g) = \frac{d}{dt}_{t = 0}g \cdot \exp(th)$. The vector fields, $X_{h}(g)$, are left-invariant and satisfy $[X_{h_{1}}, X_{h_{2}}] = X_{[h_{1}, h_{2}]}$. If $\pi:G \to G/P$ is the left coset projection, the image, $X_{h}(gP) = \pi_{\ast}(g)(X_{h})$, is a left-invariant vector field on $G/P$, which vanishes if and only if $h \in \p$ ($\p$ is the Lie algebra of $P$). If $h = r^{\alpha}e_{\alpha}$, where
\begin{align*}
&e_{\alpha} = \begin{pmatrix} 0 & 0 & 0 \\ \delta_{\alpha}\,^{j} & 0 & 0 \\ \delta_{\alpha}\,^{0} & -\omega_{\alpha q} & 0 \end{pmatrix},\,\, \text{and} \,\, x = \begin{pmatrix} 1 & 0 & 0\\ x^{j} & \delta_{i}\,^{j} & 0 \\ x^{0} & -x_{i} & 1 \end{pmatrix}P,& \\&\text{then}\,\, X_{h}(x) = \begin{pmatrix} 0 & 0 & 0 \\ r^{j} & 0 & 0 \\ r^{0} + r_{p}x^{p} & -r_{i} & 0 \end{pmatrix}.
\end{align*}
This shows that $X_{i} = X_{e_{i}} = \frac{\partial}{\partial x^{i}} + \omega_{i p}x^{p}\frac{\partial}{\partial x^{0}}$ and $X_{0} = 2X_{e_{0}} = 2\frac{\partial}{\partial x^{0}}$ constitute a left-invariant basis of $T(G/P)$. The Lie brackets are $[X_{\alpha}, X_{\beta}] = -\omega_{\alpha\beta}X_{0}$. The rank $2n-2$ left-invariant subbundle of $T(G/P)$ spanned by the vector fields $X_{i}$ is the canonical contact structure, $H$, on $G/P$. A left-invariant section of the annihilator of $H$ is given by $\theta = \tfrac{1}{2}(dx^{0} + \omega_{pq}x^{p}dx^{q})$. The forms $\theta^{i} = dx^{i}$ constitute with $\theta$ a left-invariant $\theta$-adapted coframe, and $X_{0}$ is the Reeb vector field of $\theta$. The flat model contact projective structure is induced by the Maurer-Cartan form on $G$ viewed as a Cartan connection on the bundle $G \to G/P$. The representative of the flat model contact projective structure associated to $\theta$ by Theorem \ref{theorema} is the unique $\nabla$ determined by requiring the $X_{\alpha}$ to be parallel. The ambient connection associated to this flat model contact projective structure is the usual Euclidean connection on $\rea^{2n}$.

\section{Contact Schwarzian Derivative}
\subsection{Contact Projective Schwarzian Cocycle}
There is a well defined notion of the difference tensor of two contact projective structures, $[\bar{\nabla}]$ and $[\nabla]$, as a section of the bundle $\Ct$ defined in Section \ref{notationsection}. Fix a contact one-form, $\theta$, and let $\Pi$ be the difference tensor of the representatives, $\bar{\nabla} \in [\bar{\nabla}]$ and $\nabla \in [\nabla]$, associated to $\theta$ by Theorem \ref{theorema}. Each of $\bar{\nabla}$ and $\nabla$ makes parallel $\theta$ and $\rb$, and the interior multiplication of $\rb$ in the torsion of each vanishes; these imply
\begin{equation}\label{aff0}
 \Pi_{\alpha\beta}\,^{0} = \Pi_{\alpha 0}\,^{\beta} = \Pi_{0\alpha}\,^{\beta} = 0,
\end{equation}
so that $\Pi$ may be identified with the section, $\Pi_{ij}\,^{k}$, of $\tensor^{2}(H^{\ast})\tensor H$. Let $\Tilde{\Pi}$ be the difference tensor of the representatives, $\Tilde{\bar{\nabla}} \in [\bar{\nabla}]$ and $\tilde{\nabla} \in [\nabla]$, associated to $\tilde{\theta} = f^{2}\theta$ by Theorem \ref{theorema}. Letting $\tilde{\Pi}_{\alpha\beta}\,^{\gamma}$ denote the components of $\tilde{\Pi}$ with respect to a $\tilde{\theta}$-adapted coframe and dual frame, observe that, as in \eqref{aff0}, $\tilde{\Pi}_{\alpha\beta}\,^{0} = \tilde{\Pi}_{\alpha 0}\,^{\beta} = \tilde{\Pi}_{0\alpha}\,^{\beta} = 0$. As a consequence, the components of $\tilde{\Pi}_{ij}\,^{k}$ are the same when calculated in a $\theta$-adapted coframe and dual frame as when calculated in a $\tilde{\theta}$-adapted coframe and dual frame. It is now claimed that $\tilde{\Pi}_{ij}\,^{k} = \Pi_{ij}\,^{k}$. Let $\bar{\Lambda}$ be the difference tensor of $\tilde{\bar{\nabla}}$ and $\bar{\nabla}$, let $\Lambda$ be the difference tensor of $\tilde{\nabla}$ and $\nabla$, and observe that $\tilde{\Pi}-\Pi = \bar{\Lambda} - \Lambda$. \eqref{lambdatransform} shows that $\bar{\Lambda}_{ij}\,^{k} = \Lambda_{ij}\,^{k}$, so that $\tilde{\Pi}_{ij}\,^{k} = \Pi_{ij}\,^{k}$. Hence $\Pi_{ij}\,^{k}$ is independent of the choice of $\theta$, and, consequently, it makes sense to speak of $\Pi_{ij}\,^{k}$ as the difference tensor of the contact projective structures. It can be checked that $\Pi_{[ijk]} = 0$ and $\Pi_{ij}\,^{k}$ is completely trace-free, so that $\Pi_{ij}\,^{k} \in \Ct$.

Theorem $2.2$ of \cite{Fox} describes the affine structure of the non-empty space of contact projective structures on $(M, H)$. The difference tensor, $\Pi_{ij}\,^{k}$, of two contact projective structures on $M$ admits a direct sum decomposition, $\Pi_{ij}\,^{k} = A_{ij}\,^{k} + B_{ij}\,^{k}$, where $A_{ij}\,^{k}\in \Gamma(\St)$ and $B_{ij}\,^{k} \in \Gamma(\Tt)$.
In particular, the difference of the contact torsions is $2\Pi_{[ij]}\,^{k}$, and the difference tensor of two contact projective structures with the same contact torsion satisfies $\Pi_{ijk} = \Pi_{(ijk)}$. Given $p \in M$ there is an open $U \subset M$, containing $p$, so that for any trace-free section, $\tau_{ij}\,^{k}$, defined over $U$ and satisfying $\tau_{[ijk]} = 0$, there exists in $U$ a contact projective structure with contact torsion $\tau_{ij}\,^{k}$.

The group of diffeomorphisms of $M$ acts on the space of affine connections on $M$; as a differential operator the result, $\pnabla$, of this action is given by $\pnabla = \phi^{\ast} \circ \nabla \circ (\phi^{-1})^{\ast}$. The connection $\pnabla$ is characterized by its action on vector fields, $\phi_{\ast}(\pnabla_{X}Y) = \nabla_{\phi_{\ast}(X)}\phi_{\ast}(Y)$, where the pullback of vector fields is defined by $\phi^{\ast}(X) = (\phi^{-1})_{\ast}(X)$. More generally, the definition of $\pnabla$ is made so that if $S$ is any weighted tensor on $M$, then $\pnabla \phi^{\ast}(S) = \phi^{\ast}(\nabla S)$. For any contact one-form, $\theta$, and any contactomorphism, $\phi$, there follows that $\pnabla_{(i}\theta_{j)}$ is the pullback via $\phi$ of $\nabla_{(i}\theta_{j)}$, so $\pnabla$ admits a full set of contact geodesics if and only if $\nabla$ admits a full set of contact geodesics; in this case a contact path is a contact geodesic of $\pnabla$ if and only if its image under $\phi$ is a contact geodesic of $\nabla$. Hence the group, $\G = \contacto(M, H)$, of contactomorphisms of a contact manifold acts on the space of contact projective structures by $\phi \cdot [\nabla] = [\pnabla]$. Note that if $\nabla \in [\nabla]$ is the representative associated to $\theta$ by Theorem \ref{theorema} then $\pnabla \in [\pnabla]$ is the representative associated to $\phi^{\ast}(\theta)$ by Theorem \ref{theorema}. The space of contact projective structures on $(M, H)$ is an affine space modeled on the infinite-dimensional vector space, $\ab = \Gamma(\Ct)$, and $\G$ acts on $\ab$ by $\phi \cdot \Omega = \phi^{\ast}(\Omega)$. The chain complex $C^{k}(\G; \ab)$, for $\G$ with coefficients in $\ab$, is the space of maps from $\G^{k}$ to $\ab$ with the usual coboundary of group cohomology. Precisely,
\begin{align}
&\text{for $\Omega \in C^{0}(\G, \ab) = \ab$,}& &\partial\Omega(\phi) = \phi^{\ast}(\Omega) - \Omega,\\
&\text{for $T \in C^{1}(\G, \ab)$,}& &\partial T(\psi, \phi)  = \psi^{\ast}(T(\phi)) - T(\phi \circ \psi) + T(\psi).
\end{align} 

\begin{definition}
For a contact projective structure, $(M, H, [\nabla])$, define the \textbf{contact projective Schwarzian derivative}, $\schwn \in C^{1}(\G, \ab)$, i.e. $\schwn:\G = \contacto(M, H) \to \ab = \Gamma(\Ct)$, by letting $\schwn(\phi)$ be the difference tensor of $[\pnabla]$ and $[\nabla]$. 
\end{definition}
\noindent
By definition $\schwn$ has the property that $\schwn(\phi) = 0$ if and only if $\phi$ is an automorphism of the contact projective structure represented by $[\nabla]$. 

\begin{lemma}
The contact projective Schwarzian derivative has the following properties:
\begin{enumerate}
\item (Cocycle Property). $\partial \schwn = 0$.
\item (Equivariance property). For $\Omega \in \ab$, $\schw_{[\nabla] + \Omega} - \schwn = \partial \Omega$.
\end{enumerate}
\end{lemma}
\begin{proof}
By definition $[\pnabla] - [\nabla]$ is the difference tensor of the unique representatives, $\bar{\nabla} \in [\pnabla]$ and $\nabla \in [\nabla]$, making parallel a chosen contact one-form, $\theta$. As $^{\psi}\bar{\nabla}$ is projectively equivalent to $^{\psi}(\pnabla) =\, ^{\phi \circ \psi}\nabla$, and the unique representatives making parallel $\psi^{\ast}(\theta)$ are $^{\psi}\bar{\nabla} \in [^{\psi}(\pnabla)]$ and $^{\psi}\nabla \in [^\psi\nabla]$, the difference tensor of $^{\psi}\bar{\nabla}$ and $^\psi\nabla$ is by definition $[^{\psi}(\pnabla)] - [^\psi\nabla]$, which by definition is the pullback via $\psi$ of the tensor $[\pnabla] - [\nabla]$. This shows
\begin{align}\label{equilemma}
[^{\psi}(\pnabla)] - [^\psi\nabla] = \psi^{\ast}([\pnabla] - [\nabla]) = \psi^{\ast}(\schwn(\phi)).
\end{align}
The cocycle property, $\schwn(\phi\circ \psi) = \psi^{\ast}(\schwn(\phi)) + \schwn(\psi)$, means explicitly that, for $X, Y \in \Gamma(H)$,
\begin{align*}
&\schwn(\phi \circ \psi)(X, Y) = (\psi^{-1})_{\ast}(\schwn(\phi)(\psi_{\ast}(X), \psi_{\ast}(Y))) + \schwn(\psi)(X, Y),
\end{align*}
and this follows from \eqref{equilemma}. The equivariance property follows similarly from the definitions. 
\end{proof}
\noindent
The equivariance property has as a special case the identity $\schw_{[\pnabla]} = \schwn + \partial(\schwn(\phi))$, which may be rewritten as $\schw_{[\pnabla]} = \phi \cdot \schwn$, where the natural action of $\G$ on $C^{1}(\G; \ab)$ is defined by $(\phi \cdot T)(\psi) = \phi^{\ast}(T(\phi\circ \psi\circ \phi^{-1}))$ (if $\partial T = 0$ then $\phi \cdot T = T + \partial(T(\phi))$).
\begin{proposition}
The cohomology class $[\schwn] \in H^{1}(\G, \ab)$ is non-trivial and does not depend on the choice of contact projective structure. 
\end{proposition}

\begin{proof}
The cocycle property shows that $\partial \schwn = 0$ for any $[\nabla]$, so $\schwn$ determines a class $[\schwn] \in H^{1}(\G, \ab)$. The equivariance property shows that for any other contact projective structure, $[\bar{\nabla}]$, $\schw_{[\bar{\nabla}]} - \schwn$ is a coboundary, and so the cohomology class $[\schwn]$ does not depend on the choice of $[\nabla]$. Suppose there were $\Omega \in \ab$ such that $\schwn = \partial \Omega$. Then
\begin{align}\label{interinv}
&[\pnabla] - [\nabla] = \schwn(\phi) = \partial \Omega (\phi) = \phi^{\ast}(\Omega) - \Omega.
\end{align}
Let $[\nablap] = [\nabla] -\Omega$. Then \eqref{interinv} shows $[^{\phi}\nablap] = [\pnabla]  -\phi^{\ast}(\Omega) = [\nablap]$ for all contactomorphisms $\phi$, which contradicts the finite-dimensionality of the automorphism group of the contact projective structure $[\nablap]$. (See \cite{Cap-automorphism} for a simple proof that the automorphism group of a parabolic geometry is finite dimensional). This shows that $\schwn$ is not a coboundary.
\end{proof}
The \textbf{contact Schwarzian derivative} is defined by taking $[\nabla]$ to be the flat model contact projective structure described in section \ref{flatmodelsection}. The resulting cocycle, $\schw(\phi)$, has by definition the property that $\schw(\phi) = 0$ if and only if $\phi$ is the restriction to an open neighborhood of the action on $\proj(\rea^{2n})$ of an element of $Sp(n, \rea)$. Next this $\schw(\phi)$ is expressed explicitly in terms of a local Darboux frame. Let $\nabla$ be as in Section \ref{flatmodelsection} and fix a $\theta$-adapted frame and coframe. Let $\phi$ be a contactomorphism and write $\phi^{\ast}(\theta)_{x} = c(x)\theta_{x}$. In the explicit expressions to follow it is useful to keep in mind that indices are to be interpreted with respect to $\theta$-adapted frame and coframe, and also that the Reeb field of $\phi^{\ast}(\theta)$ is $\phi^{-1}_{\ast}(\rb)$. The claim is that
\begin{align}\label{explicitschwarzian}
&\schw_{ij}\,^{k}(\phi) = \lambda_{(ij)}\,^{k} - \tfrac{1}{2n-1}(\delta_{i}\,^{k}\lambda_{(jp)}\,^{p} +\delta_{j}\,^{k}\lambda_{(ip)}\,^{p} ) = \lambda_{ij}\,^{k} - \tfrac{1}{2n}(2\delta_{(i}\,^{k}\lambda_{j)\alpha}\,^{\alpha} + \omega_{ij}\lambda^{k}\,_{\alpha}\,^{\alpha}),
\end{align}
where $\lambda$ is the difference tensor of $\pnabla$ and $\nabla$. Let $\bar{\nabla} \in [\pnabla]$ be the unique representative determined by $\theta$. Since $\pnabla \phi^{\ast}(\theta) = 0$, equation \eqref{lambdatransform} of Lemma \ref{transformationruleslemma} shows that the difference tensor, $\Xi$, of $\pnabla$ and $\bar{\nabla}$ has components $\Xi_{ij}\,^{k} = \gamma_{i}\delta_{j}\,^{k} + \gamma_{j}\delta_{i}\,^{k} + \omega_{ij}\gamma^{k}$, where $\gamma = \tfrac{1}{2}d\log{c}$. By definition, $\schw_{ij}\,^{k} = \lambda_{ij}\,^{k} - \Xi_{ij}\,^{k}$. The components $2\lambda_{[ij]}\,^{k}$ may be computed from the difference of the torsions of $\pnabla$ and $\nabla$, and this shows $\lambda_{[ij]}\,^{k} = \omega_{ij}\gamma^{k}$. From this there follows that $\schw_{ij}\,^{k} = \lambda_{ij}\,^{k} - \Xi_{ij}\,^{k} = \lambda_{(ij)}\,^{k} - \Xi_{(ij)}\,^{k}$. Because $\schw_{ij}\,^{k}$ must be completely trace free, there follows $\lambda_{(ip)}\,^{p} = (2n-1)\gamma_{i}$, and this implies the first equality of \eqref{explicitschwarzian}. Similar straightforward computations show $\lambda_{i0}\,^{0} = 2\gamma_{i}$, $\lambda_{0i}\,^{0} = 0$, and $\lambda_{i[jk]} = \gamma_{i}\omega_{jk}$, from which follows $\lambda_{i\alpha}\,^{\alpha} = \lambda_{\alpha i}\,^{\alpha} = 2n\gamma_{i}$. This justifies the second equality of \eqref{explicitschwarzian}.

The reader desirous of seeing \eqref{explicitschwarzian} expressed more explicitly should proceed as follows. Denote by $\phi^{\alpha}$ the components of $\phi$. Let $\theta^{\alpha}$, $X_{\alpha}$ be as in Section \ref{flatmodelsection}. By assumption there hold
\begin{align}\label{contactoexplicit}
&\frac{\partial \phi^{0}}{\partial x^{0}} + \omega_{pq}\phi^{p}\frac{\partial \phi^{q}}{\partial x^{0}} = c,&
&\frac{\partial \phi^{0}}{\partial x^{i}} + \omega_{pq}\phi^{p}\frac{\partial \phi^{q}}{\partial x^{i}} = c\omega_{pi}x^{p}.
\end{align}
Define $A_{\alpha}\,^{\beta}$ by $\phi_{\ast}(X_{\alpha}) = A_{\alpha}\,^{\beta}X_{\beta}\circ \phi$. Computing with \eqref{contactoexplicit} shows
\begin{align*}
&A_{i}\,^{j} = \frac{\partial \phi^{j}}{\partial x^{i}} + \omega_{ip}x^{p}\frac{\partial \phi^{j}}{\partial x^{0}},& &A_{0}\,^{0} = c,& &A_{0}\,^{i} = 2\frac{\partial \phi^{i}}{\partial x^{0}},& &A_{i}\,^{0} = 0.
\end{align*}
There holds $A_{i}\,^{p}A_{j}\,^{q}\omega_{pq} = c\omega_{ij}$, so that if $B_{i}\,^{p}A_{p}\,^{j} = \delta_{i}\,^{j}$, then $B_{i}\,^{j} = -c^{-1}A^{j}\,_{i}$, and
\begin{align*}
\pnabla_{X_{i}}X_{j} = (\phi^{-1})_{\ast}(\nabla_{\phi_{\ast}(X_{i})}\phi_{\ast}(X_{j})) = A_{i}\,^{q}X_{q}(A_{j}\,^{p})B_{p}\,^{k}X_{k}, 
\end{align*} 
so that $\lambda_{ij}\,^{k} =  A_{i}\,^{q}X_{q}(A_{j}\,^{p})B_{p}\,^{k}$ and $\lambda_{ip}\,^{p} = A_{i}\,^{q}X_{q}(A_{p}\,^{s})B_{s}\,^{p}$. 



\subsection{Relationship with Contact Hessian and Proof of Theorem \ref{maintheorem}}\label{contacthessiansection}
Given a contact projective structure, by virtue of the maximal non-integrability of the contact structure, the covariant derivatives in the directions transverse to the contact structure are completely determined by the covariant derivatives in the contact directions. Precisely, tracing the Ricci identity shows that $\nabla_{0} = \tfrac{1}{n-1}\nabla^{p}\nabla_{p} \,\mod \text{$0$th order terms}$; in particular, for $u \in \emf[\lambda]$, the Ricci identity implies $2\nabla_{[i}\nabla_{j]}u + \tau_{ij}\,^{\alpha}\nabla_{\alpha}u = 0$, and tracing this gives $\nabla^{p}\nabla_{p}u = (n-1)\nabla_{0}u$. The operator, $K:\Gamma(\emf[1]) \to \Gamma(\emfz_{ij}[1])$, defined by 
\begin{equation}\label{invop1}
K_{ij}u = \nabla_{i}\nabla_{j}u + \tfrac{1}{2}\omega_{ij}\nabla_{0}u - P_{ij}u =\nabla_{i}\nabla_{j}u + \tfrac{1}{2(n-1)}\omega_{ij}\nabla^{p}\nabla_{p}u - P_{ij}u,
\end{equation}
is an invariant differential operator; the trace free part of $\nabla_{i}\nabla_{j}u - P_{ij}u$ is $K_{ij}u$. The invariance means that $\tilde{K}_{ij}u = K_{ij}u$; this follows from the transformation rules for the covariant derivatives of $u \in \Gamma(\emf[\lambda])$, (see \cite{Fox}):
\begin{align}\label{densitytransform}
&\tilde{\nabla}_{i}u - \nabla_{i}u = \lambda\gamma_{i}u,&
&f^{2}\tilde{\nabla}_{0}u - \nabla_{0}u = \lambda\gamma_{0}u + 2\gamma^{p}\nabla_{p}u.&
\end{align}
The skew part, $K_{[ij]}u$, is first order and vanishes when the contact torsion vanishes, so it makes sense to focus attention on the \textbf{contact Hessian}, $L_{ij}u = K_{(ij)}u = \nabla_{(i}\nabla_{j)}u - P_{(ij)}u$. 

The contact Hessian is most easily derived by use of the ambient connection, $\hnabla$, of Theorem \ref{ambienttheorem}. Let $u \in \Gamma(\emf[\lambda])$, and let $\tilde{u}$ denote the corresponding homogeneity $\lambda$ function on $\form$. Using the explicit description of $\hnabla$ available in \cite{Fox} gives directly
\begin{align*}
&\hnabla_{i}d\tilde{u}_{j} = \nabla_{i}\nabla_{j}u + \tfrac{1}{2}\omega_{ij}\nabla_{0}u - \lambda P_{ij}u,&\\
&\hnabla_{i}d\tilde{u}_{0} = \nabla_{i}\nabla_{0}u -2P_{i}\,^{p}\nabla_{p}u - \lambda P_{i0}u,\\
&\hnabla_{0}d\tilde{u}_{i} = \nabla_{0}\nabla_{i}u -2P_{i}\,^{p}\nabla_{p}u -Q_{i}\,^{p}\nabla_{p}u - \lambda P_{0i}u,&\\
&\hnabla_{0}d\tilde{u}_{0} = \nabla_{0}\nabla_{0}u -2P_{0}\,^{p}\nabla_{p}u  - \lambda P_{00}u,\\
&\hnabla_{\alpha}d\tilde{u}_{\infty} = \hnabla_{\infty}d\tilde{u}_{\alpha} = (\lambda-1)\nabla_{\alpha}u,& \\
&\hnabla_{\infty}d\tilde{u}_{\infty} = \lambda(\lambda-1)u.
\end{align*}
where the explicit expressions for $P_{i0}$, $P_{0i}$, and $P_{00}$ are given in \cite{Fox}. All the right hand expressions such as $\nabla_{i}\nabla_{j}u$ should really carry a $\tilde{\,}$ which has been dropped for ease of reading. When $\tau_{ij}\,^{k} = 0$, $\hnabla$ is torsion free, so $\hnabla_{[I}d\tilde{u}_{J]} = 0$. Specializing to $\lambda = 1$, the symmetrized operator $\hnabla_{(I}d\tilde{u}_{J)}$ gives rise to the contact Hessian, and the expressions above give a second demonstration of the invariance of $K_{ij}u$. Assuming the contact torsion vanishes, and making use of the identities
\begin{align*}
&(n-1)\nabla_{0}u = \nabla^{p}\nabla_{p}u,& &n\nabla_{i}\nabla_{0}u = \nabla^{p}\nabla_{i}\nabla_{p}u + 2nP_{i}\,^{p}\nabla_{p}u,\\
&\nabla_{i}\nabla_{0}u = \nabla_{0}\nabla_{i}u,& &(n-1)\nabla_{0}\nabla_{0}u = \nabla^{p}\nabla_{0}\nabla_{p}u,
\end{align*}
gives the expressions,
\begin{align*}&\hnabla_{i}d\tilde{u}_{j} = L_{ij}u,&\\
&\hnabla_{i}d\tilde{u}_{0} = L_{i0}u = \tfrac{2}{2n-1}\nabla^{p}L_{ip}u,&\\
&\hnabla_{0}d\tilde{u}_{i} =  L_{0i}u =\tfrac{2}{2n-1}\nabla^{p}L_{ip}u,&\\ 
&\hnabla_{0}d\tilde{u}_{0} =  L_{00}u = \tfrac{2}{(n-1)(2n-1)}\nabla^{i}\nabla^{j}L_{ij}u - \tfrac{2}{n-1}P^{ij}L_{ij}u,
\end{align*}
valid when the contact torsion vanishes. This proves
\begin{lemma}\label{paralleltractorlemma}
If $(M, H, [\nabla])$ is a contact projective structure with vanishing contact torsion, then $u \in \Gamma(\emf[1])$ solves $L_{ij}u = 0$ if and only if $d\tilde{u}$ is $\hnabla$-parallel.
\end{lemma}
\noindent
The lemma shows that if $\tau_{ij}\,^{k} = 0$, solutions of $L_{ij}u = 0$ are in bijection with $\hnabla$-parallel, homogeneity $-1$ Hamiltonian vector fields on $\form$. Since $d\tilde{u}$ has homogeneity one, it can be identified with a parallel section of the cotractor bundle constructed in \cite{Fox}, and conversely any such parallel section arises in this way. Densities, $u, v \in \emf[\lambda]$ are said to be \textbf{linearly independent} if the corresponding homogeneity $\lambda$ one-forms $d\tilde{u}$ and $d\tilde{v}$ are linearly independent. Because $\hnabla$ is flat if and only if it admits locally a parallel coframe, Lemma \ref{paralleltractorlemma} implies
\begin{proposition}\label{flatcproj}
On a manifold with a contact projective structure with vanishing contact torsion there are in a neighborhood of every point $2n$ linearly independent solutions of $L_{ij}u = 0$ if and only if the contact projective structure is flat.
\end{proposition}
Locally a density, $u$, may be regarded  as a function, and with respect to the frame, $X_{\alpha}$, the local coordinate expression of the flat model contact Hessian of $u$ is 
\begin{align*}
&L_{ij}u = \frac{\partial^{2}u}{\partial x^{i} \partial x^{j}} + \omega_{ip}x^{p}\frac{\partial^{2}u}{\partial x^{j}\partial x^{0}} + \omega_{jp}x^{p}\frac{\partial^{2}u}{\partial x^{i}\partial x^{0}} + \omega_{ip}\omega_{jq}x^{p}x^{q}\frac{\partial^{2}u}{\partial x^{0}\partial x^{0}}
\end{align*}
The left action of $G$ on $G/P$ induces an embedding of Lie algebras, $\g  \to \vect(G/P)$, defined by $X^{h}(x) = \frac{d}{dt}_{t = 0}\exp(-th)\cdot x$, so that
\begin{align*}
&h = \begin{pmatrix} a & c_{q} & c_{0} \\ b^{p} & A_{q}\,^{p} & c^{p} \\ b_{0} & -b_{q} & -a \end{pmatrix}\to X^{h}(x) =\\
\begin{split}&(-\tfrac{1}{2}b_{0} + \tfrac{1}{2}A_{pq}x^{p}x^{q} + b_{p}x^{p} + ax^{0} + c_{q}x^{q}x^{0} + \tfrac{1}{2}c_{0}(x^{0})^{2})X_{0}\\
&+(-b^{p} + (a\delta_{q}\,^{p} - A_{q}\,^{p})x^{q} - c^{p}x^{0} + (c_{q}x^{q} + c_{0}x^{0})x^{p})X_{p}
\end{split}
\end{align*}
Via the following argument, which mimics an argument of M. Eastwood, \cite{Eastwood}, in the projective case, Proposition \ref{flatcproj} allows the identification of the kernel of the flat model contact Hessian with the standard representation of $\g$. The Lie derivative induces an action of the algebra of vector fields on $G/P$ on smooth sections of $\emf[\lambda]$. Restricting this action to the $X^{h}$ gives a representation of $\g$ on the space of smooth sections of $\emf[\lambda]$. Locally densities may be regarded as functions, and this action of $\g$ is represented by first order differential operators of the following form
\begin{align*}
&h \to X^{h} - \lambda(a + c_{p}x^{p} + c_{0}x^{0}) = D^{h}_{\lambda}.
\end{align*}
For any fixed $\lambda$, the map $h \to D^{h}_{\lambda}$ embeds $\g$ in the Lie algebra of differential operators. Invariantly, $D^{g}_{\lambda}$ describes the action of $\lie_{X^{g}}$ on $\emf[\lambda]$. Since the action of $G$ on $G/P$ is by automorphisms of the flat contact projective structure, it leaves $L$ invariant, and consequently, for all $g \in \g$, $L$ commutes with $\lie_{X^{g}}$ on sections of $\emf[1]$. It follows that the space, $\rep$, of smooth sections of $\emf[1]$ annihilated by $L$ is a representation space for $\g$. As the operators $X^{e_{i}}$ and $X^{e_{0}}$ are contained among the lowering operators, any lowest weight vector is constant; thus there is in $\rep$ a unique lowest weight vector up to scale. Lemma \ref{paralleltractorlemma} implies that $\rep$ is finite-dimensional, so the theorem of the highest weight implies that $\rep$ is irreducible, and straightforward explicit computation 
shows that $\rep$ must be of lowest weight $-1$, so the standard representation. It follows that $\rep$ is generated by applying to a lowest weight vector the raising operators, and it is straightforward to check that one obtains in this way elements of the form $a + a_{p}x^{p} + a_{0}x^{0}$, for some constants, $a, a_{p}, a_{0}$. Precisely, applying to $1$ any $D^{g}_{1}$ gives an element having such a form, and applying to this element any $D^{g'}_{1}$ gives another element of the same form. 
This completes the proof of the following corollary of Proposition \ref{flatcproj}.
\begin{corollary} \label{repspacelemma}
For the flat model contact projective structure, the space of smooth solutions of $L_{ij}u = 0$, is the standard $2n$-dimensional representation of $Sp(n, \rea)$. Precisely, the solutions of $L_{ij}u = 0$ all have locally the form $a + a_{p}x^{p} + a_{0}x^{0}$.
\end{corollary}
The PDE \eqref{schwarzianpde} appearing after Theorem \ref{maintheorem} is simply the explicit expression of the contact Hessian with respect to a $\nabla$-parallel, $\theta$-adapted coframe. Corollary \ref{repspacelemma} implies that there is a basis of solutions of \eqref{schwarzianpde} constituted by $2n$ linearly independent densities obtained as the pullbacks via $\phi$ of the densities represented in local coordinates by the functions $1$ and $x^{\alpha}$. The ratio of two non-vanishing densities is a function, and the ratios of the pullbacks, $\frac{\phi^{\ast}(x^{\alpha})}{\phi^{\ast}(1)}$ are the coordinate functions, $\phi^{\alpha}$, of $\phi$. This shows that there is a basis of solutions, $f^{\infty}, f^{\alpha}$, of $L_{ij}u = 0$, so that $\phi^{\alpha} = \frac{f^{\alpha}}{f^{\infty}}$. Consequently knowing the Schwarzian derivative of a contactomorphism is enough to recover locally the contactomorphism by solving a system of PDEs. This completes the proof of Theorem \ref{maintheorem}.

\subsection{Integrability Conditions}\label{integrabilitysection}
By solving a linear system of PDE there may be constructed locally from any section $\Pi \in \Gamma(\St)$ satisfying an appropriate integrability condition, a contactomorphism, $\phi$, such that $\schw(\phi) = \Pi$. This is best understood as providing a means of constructing explicitly a developing map for a flat contact projective structure. In this section let $\nabla$ denote the representative of the standard flat contact projective structure associated to $\theta = dx^{0} + \omega_{pq}x^{p}dx^{q}$ and described in Section \ref{contacthessiansection}. By construction $\schw(\phi)$ satisfies some integrability conditions implied by the vanishing of curvatures. If the dimension is at least five, the flatness of $[\pnabla]$ is equivalent to $\tau_{ij}\,^{k} = 0 = W_{ijk}\,^{l}$. The vanishing of the contact torsion means $\schw_{(ijk)}(\phi) = \schw_{ijk}(\phi)$, and the vanishing of the contact projective Weyl tensor gives a PDE satisfied by $\schw(\phi)$. On the other hand, given any $A_{ij}\,^{k} \in \Gamma(\St)$ satisfying these integrability conditions, define a contact projective structure by requiring that its difference tensor with the flat contact projective structure be $A_{ij}\,^{k}$. The resulting contact projective structure is flat, so is locally equivalent to the flat projective structure, and so there must be defined locally a contactomorphism $\phi$ such that $A_{ij}\,^{k} = \schw_{ij}\,^{k}(\phi)$. Moreover, by the discussion following Corollary \ref{repspacelemma}, $\phi$ may be constructed explicitly by defining $\phi^{\alpha} = \frac{f^{\alpha}}{f^{\infty}}$ where $f^{\infty}, f^{\alpha}$ are a basis of solutions of $L_{ij}u = 0$, $L$ being the contact Hessian of the new contact projective structure. In three dimensions the integrability condition on a tensor to be the Schwarzian derivative of some contactomorphism comes instead from the PDE given by the vanishing of the tensor $C_{ijk}$ associated to $[\pnabla]$. Everything else works exactly as in the higher dimensional case. The resulting integrability conditions can be written down explicitly, and to convince the reader that such explicit expressions are rather useless, they are recorded here. Precisely, if $2n - 1 > 3$, $A_{ij}\,^{k} \in \Gamma(\St)$ is the contact Schwarzian of some locally defined contactomorphism, $\phi$, if and only if 
$A_{ij}\,^{k}$ satisfies the integrability condition
\begin{align*}&
 0 = 2n(\nabla_{[i}A_{j]kl} - A_{k[i}\,^{p}A_{j]lp})  - \omega_{ij}(\nabla_{p}A_{kl}\,^{p} - A_{kq}\,^{p}A_{pl}\,^{q})\\
&\notag + \nabla_{p}A_{k[i}\,^{p}\omega_{j]l} -A_{kq}\,^{p}A_{p[i}\,^{q}\omega_{j]l}+ \nabla_{p}A_{l[i}\,^{p}\omega_{j]k} -A_{lq}\,^{p}A_{p[i}\,^{q}\omega_{j]k},
\end{align*}
which amounts to requiring that $\nabla_{[i}A_{j]k}\,^{l} - A_{k[i}\,^{p}A_{j]p}\,^{l}$ be completely trace free. The condition in dimension $3$ is similar, though more involved, so is omitted. The remainder of this section describes a reformulation of these integrability conditions.

Theorem C of \cite{Fox} associates to each contact manifold, $(M, H)$, a $P$ principal bundle $\pi:\bun \to M$, the bundle of filtered projective symplectic frames in the tractor bundle, and to each contact projective structure on $M$ a canonical regular $(\g, P)$ Cartan connection, $\eta$, on $\bun$, such that $\eta$ is normal if and only if the contact projective structure is contact torsion free. The total space of the bundle $\form$ is recovered as a quotient $\bun/\tilde{P}$, where $\tilde{P} \subset P$ is the subgroup preserving a fixed vector in the standard representation, $\standrep$, of $G$. Moreover, $\eta$ is normalized by the requirement that the covariant differentiation induced by $\eta$ on the associated bundle $\bun \times_{\tilde{P}} \standrep \simeq T\form \to \form$ should be the ambient connection of Theorem \ref{ambienttheorem}. 

Let $\hpnabla$ and $\hnabla$ be the ambient connections associated to $[\pnabla]$ and $[\nabla]$, and let $\hschwn$ be their difference tensor. The components of $\hschwn$ are easily computed:
\begin{align*}
&\hschwn_{ij}\,^{k} = \schwn_{ij}\,^{k},& &\hschwn_{\alpha \beta}\,^{\infty} = \bar{P}_{\alpha \beta} - P_{\alpha \beta},&  \\ &\hschwn_{IJ}\,^{0} = 0,& &\hschwn_{0j}\,^{k} = 2\bar{P}_{j}\,^{k} - 2P_{j}\,^{k} + \bar{Q}_{j}\,^{k} - Q_{j}\,^{k},&\\
  &\hschwn_{\infty I}\,^{J} = 0 = \hschwn_{I\infty}\,^{J},&&\hschwn_{j0}\,^{k} = 2\bar{P}_{j}\,^{k} - 2P_{j}\,^{k},& 
\end{align*}
where the explicit expressions of the differences $\bar{P}_{\alpha\beta} - P_{\alpha\beta}$ and $\bar{Q}_{\alpha\beta} - Q_{\alpha\beta}$ in terms of $\schwn_{ij}\,^{k}$ and its covariant derivatives are omitted to save space. If $\bar{\eta}$ and $\eta$ are any $(\g, P)$ Cartan connections on $\bun$, their difference, $S = \bar{\eta} - \eta$, is a horizontal, $P$-equivariant, $\g$-valued one-form on $\bun$, so descends to a $\g$-valued one-form on $M$. If $\bar{\eta}$ and $\eta$ are the canonical Cartan connections associated to contact projective structures, $[\bar{\nabla}]$ and $[\nabla]$, on $M$, then the components of their difference $S = \bar{\eta} - \eta$ may be identified straightforwardly (via the tractor formalism described in \cite{Fox}) with the components of $\hschwn$ computed above. In general if $\Omega = d\eta + \eta \wedge \eta$ is the curvature of $\eta$, then $\bar{\Omega} - \Omega = \dn S + S \wedge S$, where $d^{\eta}$ is the twisted exterior derivative defined on  $\g$-valued $k$-forms by $\dn \Psi = d\Psi + \eta \wedge \Psi + (-1)^{k+1}\Psi \wedge \eta$. This $\dn$ has the nice property that $(\dn)^{2}(\Psi) = \Omega \wedge \Psi - \Psi \wedge \Omega$, so that if $[\bar{\nabla}]$ and $[\nabla]$ are flat then $\dn S + S\wedge S = 0$. As is easily checked by rewriting it in terms of the components of $\hschw$, the equation $\dn S + S\wedge S = 0$ encapsulates the integrability conditions imposed on $\schw$. Conversely, given the Cartan connection, $\eta$, associated to a flat contact projective structure, and given a $\g$-valued one-form, $S$, on $M$ satisfying $\dn S + S\wedge S = 0$, the Cartan connection $\bar{\eta} = \eta + \bar{S}$ (where $\bar{S}$ is the horizontal lift of $S$), will be also flat, and consequently there may be constructed locally a contactomorphism, $\phi$, so that $\bar{\eta}$ is induced from $\eta$ by pullback via the principal bundle automorphism of $\G$ induced by $\phi$. Working at the level of the ambient connections shows that $S$ is completely and explicitly determined by $\schw(\phi)$, though again there seems no point in writing out explicitly the identifications.

This point of view hints crudely at a connection with generalized BGG sequences. Section $8$ of \cite{Calderbank-Diemer} contains a relevant discussion of moduli of flat parabolic geometries.

\subsection{Relationship with Contact Path Geometries}\label{contactpathsection}
In \cite{Sato} and \cite{Ozawa-Sato} the contact Schwarzian is derived in the context of the equivalence of three-dimensional contact path geometries; in this section contact projective structures will be discussed from the point of view of contact path geometries in order to elucidate the relationship between the formulas of this paper and the formulas of \cite{Ozawa-Sato}. A complete discussion of contact path geometries would take much more space, so is deferred to another place. Here the basic facts are stated mostly without proof.

Each contact path admits a canonical lift to a one-dimensional submanifold of the total space of the projectivized contact distribution, $\pi:\proj(H) \to M$, and the lifts of all the contact paths in a contact path geometry foliate $\proj(H)$. The prolongation of $H$ is the tautological bundle $E\to \proj(H)$ defined by $E_{L} = \pi_{\ast}(L)^{-1}(L)$. The leaves of the foliation determined by a contact path geometry are tangent to $E$ and transverse to the vertical subbundle $V = \ker \pi_{\ast}$. A bundle automorphism of $\proj(H)$ preserves $E$ if and only if it is the lift of a contactomorphism of $(M, H)$. Consequently a contact path geometry may be reformulated as a splitting, $E = V \oplus W$, the foliation comprising the integral manifolds of $W$. For example, every three-dimensional contact path geometry is locally equivalent to the contact path geometry determined by the solutions of a third order ODE considered modulo contact transformations; in \cite{Chern}, S.-S. Chern solved the local equivalence problem for these structures. 

Let $\theta = \tfrac{1}{2}(dz + x^{\infty}dx^{0} - x^{0}dx^{\infty} + \omega_{pq}x^{p}dx^{q})$. A frame spanning $H$ is given by $X_{\infty} = \frac{\partial}{\partial x^{\infty}} + x^{0}\frac{\partial}{\partial z}$, $X_{0} = \frac{\partial}{\partial x^{0}} - x^{\infty}\frac{\partial}{\partial z}$, $X_{i} = \frac{\partial}{\partial x^{i}} + \omega_{ip}x^{p}\frac{\partial}{\partial z}$. The Reeb field is $\rb = \frac{\partial}{\partial z}$. Let capital Latin indices run over $\{\infty, 1, \dots, 2n, 0\}$. (Note that the $0$ index no longer corresponds to the Reeb direction). The representative of the flat model contact projective structure associated to $\theta$ is the unique affine connection, $\nabla$, defined by requiring the left-invariant frame, $X_{\infty}, X_{i}, X_{0}, \rb$, to be parallel. The one-forms $dx^{I}$ span $H^{\ast}$ and coordinates on the fibers of $H$ are defined by $a^{I}X_{I} = dx^{I}(X)$. Define, by $u^{\alpha} = \tfrac{a^{\alpha}}{a^{\infty}}$, coordinates in a chart on the fibers of $\proj(H)$ on which $a^{\infty} \neq 0$. Write $u_{i} = u^{p}\omega_{pi}$. A frame in $T\proj(H)$ is given by:
\begin{align*}
&A_{i} = \frac{\partial}{\partial u^{i}} - u_{i}\frac{\partial}{\partial u^{0}},& &T_{-1, 0} = X_{\infty} + u^{\alpha}X_{\alpha},& & E_{i} = X_{i} - u_{i}X_{0},&\\
 &T_{0, -2} = \frac{\partial}{\partial u^{0}},&
 &T_{-1, -2} = X_{0},& &T_{-2, -2} = \frac{\partial}{\partial z}.
\end{align*} 
(in three dimensions the $A_{i}$ and $E_{i}$ should be omitted and the other vector fields should be relabeled). Each fiber $\proj(H_{L})$ is the projectivization of a symplectic vector space, so has a canonical contact structure, and this determines a rank $2n-4$ subbundle, $U \subset V$, spanned by the vector fields $A_{i}$. With $T_{0,-2}$, the $A_{i}$ span $V$, and with also $T_{-1, 0}$ they span $E$. The most general vector field spanning $W$ has the form 
\begin{align}\label{cpathform}
X = C(T_{-1, 0} + f^{0}T_{0, -2} + f^{p}A_{p}),
\end{align} 
for some non-vanishing function $C$, and it is usually convenient to choose $X$ so that $C = 1$ (here $f_{p} = f^{q}\omega_{qp}$). There is a second tautological bundle, $E^{\perp}\to \proj(H)$, defined by $E^{\perp}_{L} = \pi_{\ast}(L)^{-1}(L^{\perp})$, where $L^{\perp}$ denotes the skew complement in $H$ of $L$ (with respect to the conformal symplectic structure on $H$). Some motivation for the following definition will be provided by Lemma \ref{contactpathprojlemma}. 
\begin{definition}\label{vancontor}
A contact path geometry has \textbf{vanishing contact torsion} if for any choice of $X$ spanning $W$ and any section $A$ of $U$, the iterated bracket, $[X, [X, A]]$ is contained in $E^{\perp}$.
\end{definition}
\noindent
Since $X$ and $[X, A]$ are sections of $E^{\perp}$, the given condition does not depend on the choices of $X$ and $A$. Since for three-dimensional contact path geometries, $E^{\perp} = E$ and $U$ is trivial, three-dimensional contact path geometries necessarily have vanishing contact torsion. 

\begin{remark}\label{torsionidentityremark}
For the particular choices of $A_{i}$ and $X$ given above, explicit computation utilizing the fact that $E^{\perp}$ is spanned by $A_{i}$, $T_{0, -2}$, $X$, and $E_{i}$ shows that the condition of Definition \ref{vancontor} occurs if and only if $3f_{p} + A_{p}(f^{0}) = 0$. 
\end{remark}

Fix a contact projective structure, $[\nabla]$, and let $\nabla$ be the representative associated to $\theta$. Let $\Gamma_{IJ}\,^{K}(x^{I}, z)$ be the Christoffel symbols of $\nabla$ in the contact directions. Note that $\tau_{IJ}\,^{K} = 2\Gamma_{[IJ]}\,^{K}$ and that $\Gamma_{IJK} = \Gamma_{I(JK)}$ (because $\nabla\omega = 0$). The coefficients $\Gamma_{IJK}$ are completely determined by $\tau_{IJK}$ and $\Gamma_{(IJK)}$ by $\Gamma_{IJK} = \Gamma_{(IJK)} + \tfrac{1}{2}\tau_{IJK} - \tfrac{1}{6}\tau_{KIJ} - \tfrac{1}{6}\tau_{KJI}$. In particular $[\nabla]$ has vanishing contact torsion if and only if $\Gamma_{(IJK)} = \Gamma_{IJK}$. The equations of the contact geodesics of $[\nabla]$ are
\begin{align}
&\label{cgeodinf}\ddot{x}^{\infty} + \Gamma_{(\alpha\beta)}\,^{\infty}\dot{x}^{\alpha}\dot{x}^{\beta} + 2\Gamma_{(\alpha \infty)}\,^{\infty}\dot{x}^{\alpha}\dot{x}\,^{\infty} + \Gamma_{\infty\infty}\,^{\infty}(\dot{x}^{\infty})^{2} = 0,\\
&\label{cgeodgamma}\ddot{x}^{\gamma} + \Gamma_{(\alpha\beta)}\,^{\gamma}\dot{x}^{\alpha}\dot{x}^{\beta} + 2\Gamma_{(\alpha \infty)}\,^{\gamma}\dot{x}^{\alpha}\dot{x}\,^{\infty} + \Gamma_{\infty\infty}\,^{\gamma}(\dot{x}^{\infty})^{2} = 0,
\end{align}
subject to the non-holonomic constraint $\dot{z} + x^{\infty}\dot{x}^{0} - x^{0}\dot{x}^{\infty} + \omega_{pq}x^{p}\dot{x}^{q} = 0$. Eliminating the variable $t$ from \eqref{cgeodinf} and \eqref{cgeodgamma} by taking $x^{\infty}$ as the independent variable; writing $\dot{x}^{\alpha} = \frac{d x^{\alpha}}{d x^{\infty}}$ and $\dot{z} = \frac{dz }{dx^{\infty}}$; and substituting \eqref{cgeodinf} into \eqref{cgeodgamma} gives the system of ordinary differential equations,
\begin{align}
&\label{cpathcproj} \ddot{x}^{\gamma} = \dot{x}^{\alpha}\dot{x}^{\beta}\dot{x}^{\gamma}\Gamma_{(\alpha\beta)}\,^{\infty}+ (2\delta_{\beta}\,^{\gamma}\Gamma_{(\alpha\infty)}\,^{\infty} - \Gamma_{(\alpha\beta)}\,^{\gamma}) \dot{x}^{\alpha}\dot{x}^{\beta} \\ &\notag  + (\delta_{\alpha}\,^{\gamma}\Gamma_{\infty\infty}\,^{\infty} - 2\Gamma_{(\alpha\infty)}\,^{\gamma})\dot{x}^{\beta} - \Gamma_{\infty\infty}\,^{\gamma},\\
&\notag 0 =  \dot{z} + x^{\infty}\dot{x}^{0} - x^{0} + \omega_{pq}x^{p}\dot{x}^{q}.
\end{align}
The corresponding subbundle $W$ is spanned by the vector field
\begin{align*}
&X = T_{-1, 0} + \left(\Gamma_{(\alpha\beta)}\,^{\infty}u^{\alpha}u^{\beta}u^{\gamma} + (2\delta_{\beta}\,^{\gamma}\Gamma_{(\alpha\infty)}\,^{\infty} - \Gamma_{(\alpha\beta)}\,^{\gamma})u^{\alpha}u^{\beta}\right)\tfrac{\partial}{\partial u^{\gamma}}\\
&\notag +\left((\delta_{\alpha}\,^{\gamma}\Gamma_{\infty\infty}\,^{\infty} - 2\Gamma_{(\alpha\infty)}\,^{\gamma})u^{\alpha} - \Gamma_{\infty\infty}\,^{\gamma} \right)\tfrac{\partial}{\partial u^{\gamma}}
\end{align*}
From \eqref{cpathform}, $\Gamma_{IJ}\,^{\infty} = \Gamma_{IJ0}$, and $\Gamma_{IJ}\,^{0} = - \Gamma_{IJ\infty}$, there follow
\begin{align*}
f_{p} =& u_{p}(u^{\alpha}u^{\beta}\Gamma_{(\alpha\beta)0} + 2u^{\alpha}\Gamma_{(\alpha\infty)0} + \Gamma_{\infty\infty 0}) \\ &- (u^{\alpha}u^{\beta}\Gamma_{(\alpha\beta)p} + 2u^{\alpha}\Gamma_{(\alpha\infty)p} + \Gamma_{\infty\infty p}),\\
f^{0} + u^{p}f_{p} =& u^{\alpha}u^{\beta}(u^{0}\Gamma_{(\alpha\beta)0}+ \Gamma_{(\alpha\beta)\infty})  + 2u^{\alpha}(u^{0}\Gamma_{(\alpha\infty)0}+ \Gamma_{(\alpha\infty)\infty}) \\& + (u^{0}\Gamma_{\infty\infty 0} + \Gamma_{\infty\infty\infty}),
\end{align*}
from which follows (after a bit of manipulation)
\begin{align}\label{f0form}
&f^{0} = u^{\alpha}u^{\beta}u^{\gamma}\Gamma_{(\alpha\beta\gamma)}+ 3u^{\alpha}u^{\beta}\Gamma_{(\alpha\beta\infty)} + 3u^{\alpha}\Gamma_{(\alpha\infty\infty)} + \Gamma_{\infty\infty\infty}.
\end{align}
It is straightforward to check that
\begin{align*}
&3f_{p} + A_{p}(f^{0}) = \\ &u^{\alpha}u^{\beta}(\tau_{p(\alpha\beta)} - u_{p}\tau_{0(\alpha\beta)}) + u^{\alpha}(\tau_{p(\alpha\infty)} - u_{p}\tau_{0(\alpha\infty)}) + (\tau_{p\infty\infty} - u_{p}\tau_{0\infty\infty}),
\end{align*}
Differentiating this repeatedly in the $u^{\alpha}$ variables shows that the components $\tau_{I(JK)}$ are explicitly expressible in terms of the functions $3f_{p} + A_{p}(f^{0})$ and their derivatives in the $u^{\alpha}$ variables of order not more than three. It follows that $3f_{p} + A_{p}(f^{0}) = 0$ if and only if $0 = \tau_{I(JK)}$. The $Sp(n-1, \rea)$ representation space of trace-free tensors with symmetries $a_{[ijk]} = 0$, $a_{[ij]k} = a_{ijk}$ is canonically isomorphic to the space of trace-free tensors with symmetries $b_{(ijk)} = 0$ and $b_{i(jk)} = b_{ijk}$, (see \cite{Weyl}), and this implies that $\tau_{I(JK)} = 0$ if and only if $\tau_{IJK} = 0$. By Remark \ref{torsionidentityremark} this proves
\begin{lemma}\label{contacttorsions}
The contact path geometry induced by a contact projective structure has vanishing contact torsion in the sense of Definition \ref{vancontor} if and only if it has vanishing contact torsion as a contact projective structure. 
\end{lemma}

\begin{proposition}\label{contactpathprojlemma}
Write $f_{\alpha_{1}\dots\alpha_{k}} = \frac{\partial^{k}f}{\partial u^{\alpha_{1}}\dots \partial u^{\alpha_{k}}}$ for a smooth function $f(t, x^{\alpha}, z, u^{\alpha})$. A contact torsion free contact path geometry is a contact projective structure if and only if the contact paths are representable locally as the solution curves of a system of ordinary differential equations of the form:
\begin{align*}
&\ddot{x}^{p}  = \tfrac{1}{3}(\dot{x}^{p}f_{0}-\omega^{pq}f_{q}),& &\ddot{x}^{0} = f -\tfrac{1}{3}\dot{x}^{p}f_{p},& &\dot{z} + t\dot{x}^{0} - x^{0} + \omega_{pq}x^{p}\dot{x}^{q} = 0,
\end{align*}
where $f(t, x^{\alpha}, z, u^{\alpha})$ satisfies $f_{\alpha\beta\gamma\sigma} = 0$. In this case the Christoffel symbols of the contact projective structure may be recovered from the function $f$ by the following formulas:
\begin{align}
&\label{fabc} \Gamma_{(\alpha\beta\gamma)} = \frac{1}{6}f_{\alpha\beta\gamma}, \qquad \Gamma_{(\alpha\beta\infty)} = \frac{1}{6}f_{\alpha\beta}  - \frac{1}{6}u^{\gamma}f_{\alpha\beta\gamma},&
\\&\Gamma_{(\alpha\infty\infty)} =  \frac{1}{3}f_{\alpha}  - \frac{1}{3}u^{\beta}f_{\alpha\beta} + \frac{1}{6}u^{\beta}u^{\gamma}f_{\alpha\beta\gamma},
\\&\label{fiii} \Gamma_{\infty\infty\infty} = f - u^{\alpha}f_{\alpha} + \frac{1}{2}u^{\alpha}u^{\beta}f_{\alpha\beta} - \frac{1}{6}u^{\alpha}u^{\beta}u^{\gamma}f_{\alpha\beta\gamma}.
\end{align}
A three dimensional contact path geometry is a contact projective structure if and only if the contact paths are representable locally as the solution curves of either the system of ordinary differential equations:
\begin{align}
&\label{cp3norm1}\frac{d^{2}y}{dx^{2}} = f(x, y, z, \frac{dy}{dx}),&
&\frac{dz}{dx} = y -x\frac{dy}{dx} 
\end{align}
where $f$ is a cubic polynomial in $\frac{dy}{dx}$ and the contact form is $\theta = dz + xdy - ydx$; or of $\frac{d^{3}z}{dq^{3}} = g(q, z, \frac{dz}{dq}, \frac{d^{2}z}{dq^{2}})$, where $g$ is a cubic polynomial in $\frac{d^{2}z}{dq^{2}}$, $p = \frac{dz}{dq}$, and the contact form is $\theta = dz - pdq$.
\end{proposition}

\begin{proof}
\eqref{f0form} and the proof of Lemma \ref{contacttorsions} show that a contact torsion free contact projective structure is completely determined by $f$ satisfying $f_{\alpha\beta\gamma\sigma} = 0$. Differentiating \eqref{f0form} repeatedly recovers the $\Gamma_{(IJK)}$ from $f$ as in \eqref{fabc}-\eqref{fiii}. In the three-dimensional case, $\Gamma_{IJK} = \Gamma_{(IJK)}$ and the contact path geometry is completely determined by $f = f^{0}$; the first of the equations \eqref{cpathcproj} becomes
\begin{align}
&\label{3dcpflat} \frac{d^{2}x^{0}}{d(x^{\infty})^{2}} = (\frac{dx^{0}}{dx^{\infty}})^{3}\Gamma_{000} + 3(\frac{dx^{0}}{dx^{\infty}})^{2}\Gamma_{00\infty} + 3(\frac{dx^{0}}{dx^{\infty}})\Gamma_{\infty\infty 0} + \Gamma_{\infty\infty\infty}.
\end{align}
The final claim in the three-dimensional case is proved by analogous local coordinate computations in a frame suitably adapted to the specified contact one-form.
\end{proof}

\noindent
Proposition \ref{contactpathprojlemma} can be related to the contact Schwarzian derivative in the following manner. The flat model contact path geometry in three-dimensions is given by the graphs, in the space with variables $(t, x, z)$ and contact one-form $\theta = \tfrac{1}{2}(dz + tdx - xdt)$, of the solution curves of the system of equations, $\frac{d^{2}x}{dt^{2}} =0$, and $\frac{dz}{dt} + t\frac{dx}{dt} - x= 0$, which is the three-parameter family of lines of the form $(t, at + b, bt + c)$. If $\phi$ is a contactomorphism, then the images under $\phi^{-1}$ of the contact lines of the flat model contact projective structure are the contact geodesics of the contact projective structure $[\pnabla]$. Lemma \ref{contactpathprojlemma} and the definition of the contact Schwarzian, $\schw(\phi)$, show that these contact geodesics satisfy the equations \eqref{3dcpflat} with $\schw_{IJK}$ replacing $\Gamma_{IJK}$,
where the components of $\schw(\phi)$ are written with respect to the frame $X_{\infty}, X_{0}, \rb$. If the coordinates $(p, q, z)$ and the contact form $ dz - pdq$ were used instead for the local coordinate expressions, there would be obtained instead a single third order ODE cubic in the variables $\frac{d^{2}z}{dq^{2}}$, and the coefficients of this cubic polynomial would be constant multiples of the components of $\schw(\phi)$. These components of $\schw(\phi)$ would be exactly the four functions $P, Q, R, S$ said in \cite{Ozawa-Sato} to constitute the Schwarzian derivative of $\phi$.

In dimension at least five, \eqref{f0form} shows that a contactomorphism, $\phi$, determines locally on the total space of $\proj(H)$, a function 
\begin{align*}
&f = u^{\alpha}u^{\beta}u^{\gamma}\schw_{(\alpha\beta\gamma)}+ 3u^{\alpha}u^{\beta}\schw_{(\alpha\beta\infty)} + 3u^{\alpha}\schw_{(\alpha\infty\infty)} + \schw_{\infty\infty\infty}.
\end{align*}
Equations \eqref{fabc}-\eqref{fiii} of Proposition \ref{contactpathprojlemma} show that any such $f$ determines a contact path geometry, and if this $f$ satisfies some non-trivial integrability condition (following from the integrability condition satisfied by $\schw(\phi)$), then the resulting contact path geometry arises as the pullback via some explicitly constructible contactomorphism $\phi$ of the flat model contact path geometry.

\end{document}